\documentclass[a4paper, 12pt]{article}

\catcode`\ç\active\def ç{{\c c}} \catcode`\Ç\active\def Ç{{\c C}}
\catcode`\ð\active\def ð{{\u g}} \catcode`\Ð\active\def Ð{{\u G}}
\catcode`\ý\active\def ý{{\i}} \catcode`\Ý\active\def Ý{{\. I}}
\catcode`\ö\active\def ö{{\" o}} \catcode`\Ö\active\def Ö{{\" O}}
\catcode`\þ\active\def þ{{\c s}} \catcode`\Þ\active\def Þ{{\c S}}
\catcode`\ü\active\def ü{{\" u}} \catcode`\Ü\active\def Ü{{\" U}}

\begin{document}
\setlength{\baselineskip}{1.3\baselineskip}
\begin{center}
\textbf{CLAIMS ON PRIMORIAL PRIME NUMBERS} \\[0.25in]
TÜRKER ÖZSARI \\[0.15in]
Department of Mathematics, \\[0pt]
\smallskip Koç University, Sarýyer, Ýstanbul, 34450 \\[0pt]
\smallskip tozsari@ku.edu.tr \\[0.15in]
\today
\\[0.5in]
\textbf{Abstract}
\end{center}
\noindent {I give some claims on primorial prime numbers for interested readers in number theory.}
\section{Introduction}
\noindent \textbf{Definition 1.1} A prime number $p$ is said to be \textbf{primorial} if it is one more or less than the product of the first $n$ prime numbers for some $n\in \{1,2,...\}$.
\noindent \textbf{Definition 1.2} A prime number $p_+$ is said to be \textbf{plus-primorial} if it is one more than the product of the first $n$ prime numbers for some $n\in \{1,2,...\}$.
\noindent \textbf{Definition 1.3} A prime number $p_-$ is said to be \textbf{minus-primorial} if it is one less than the product of the first $n$ prime numbers for some $n\in \{1,2,...\}$.
\noindent \textbf{Example 1.4} 

\noindent 3 = 1 + 2 is plus-primorial.

\noindent 5 = -1 + 2.3 minus-primorial.

\noindent 7 = 1 + 2.3 is plus-primorial.

\noindent 29 = -1 + 2.3.5 is minus-primorial.

\noindent 31 = 1 + 2.3.5 is plus-primorial.

\noindent 211 = 1 + 2.3.5.7 is plus-primorial.

\noindent \textbf{Remark 1.5} An arbitrary number which is one more or less than the product of the first $n$ prime numbers for some $n\in \{1,2,...\}$ can be either a plus-primorial, a minus-primorial or a composite number.

\noindent \textbf{Example 1.6} 

\noindent 211 = 1 + 2.3.5.7 is a plus-primorial prime number.

\noindent 30031 = 59.509 = 1 + 2.3.5.7.11.13 is a composite number.

\noindent 209 = 11.19 = -1 + 2.3.5.7 is a composite number.

\noindent \textbf{Proposition 1.7} A primorial can not be both plus and minus primorial.

\noindent \textbf{Proof} : Let $p$ be both a plus and minus-primorial prime number.  Then p can be written in these two forms: $p = 1+p_1p_2...p_k = -1+p_1p_2...p_n$ for some $k$ and $n$ such that $n>k$.  Then $2 = p_1p_2...p_k(p_{k+1}...p_n-1)$.  However, right hand side is certainly greater than 2 which contradicts that $p$ is both plus and minus-primorial prime number.  Q.E.D.  

\section{My Claims}

\noindent \textbf{Claim 2.1} There are infinitely many plus-primorial prime numbers.

\noindent \textbf{Claim 2.2} There are infinitely many minus-primorial prime numbers.

\noindent \textbf{Claim 2.3} There are infinitely many prime numbers which are neither plus nor minus-primorial.

\noindent \textbf{Claim 2.4} Between any two plus-primorial prime numbers, there exists at least one prime number.

\noindent \textbf{Claim 2.5} Between any two minus-primorial prime numbers, there exists at least one prime number which is not plus-primorial.

\end{document}